\documentclass[10pt,draft,reqno]{amsart}
     \makeatletter
     \def\section{\@startsection{section}{1}%
     \z@{.7\linespacing\@plus\linespacing}{.5\linespacing}%
     {\bfseries
     \centering
     }}
     \def\@secnumfont{\bfseries}
     \makeatother
\newtheorem{theorem}{Theorem}[section]

\theoremstyle{definition}

\theoremstyle{remark}

\numberwithin{equation}{section}
\DeclareMathOperator{\dif}{\textnormal{d}\!}

\begin{document}

\title[Three dimensional polyelectrolytes]{On a remark of de Gennes 
concerning three-dimensional polyelectrolytes}

\author{Shiquan Li}
\address{Shiquan Li: Department of Mathematics, University of Rochester, 
Rochester, NY 14627, USA}
\email{shiquanli@rochester.edu}
\urladdr{https://shiquan-li.github.io}

\author{Carl Mueller}
\address{Carl Mueller: Department of Mathematics, University of Rochester, 
Rochester, NY 14627, USA}
\email{carl.e.mueller@rochester.edu}
\urladdr{https://sites.google.com/view/carlmuellerspage/home}



\begin{abstract}
This work is inspired by a remark of de Gennes \cite{deGennes79} about 
polyelectrolytes, which are charged polymers.  A common model for a 
polymer is a self-avoiding or self-repelling random walk or Brownian 
motion.  For polyelectrolytes, the repelling potential is the Coulomb 
potential arising from pairs of charged particles.  We show that in the 
continuous case of a Brownian motion in three dimensions, the spread of 
the polymer, in particular the radius of gyration of a polyelectrolyte
of length $T$ grows linearly with $T$, up to logarithmic corrections.  
\end{abstract}
\maketitle

\section{Introduction}
\label{section:introduction}

Self-avoiding random walks have been intensively studied by mathematicians,
physicists, and chemists, see \cite{MS13,dH09,BDGS12}.  In particular, 
such walks are used to model polymers. One of the most widely studied 
properties is the spread of the polymer, measured by the end-to-end 
distance, and also by the radius of gyration, see \eqref{eq:radius}.    A 
simple random walk up to time $N$ has a spread proportional to $\sqrt{N}$, 
but for self-avoiding walks, the spread is typically of order $N^\nu$, 
with $\nu>1/2$.  A well-known hard problem for self-avoiding walks is to 
find the value of $\nu$ in dimensions $d\in\{2,3\}$.  

A standard reference for the theory of polymers is de Gennes 
\cite{deGennes79}, which gives a broad overview of the subject.  One of 
the lesser-known models in \cite{deGennes79} deals with polymer chains of 
charged particles, known as polyelectrolytes, see \cite{deGennes79} 
Section XI.2.1, page 299 and also \cite{PVG77}.  In such models, the 
potential between particles is long-range, whereas more commonly studied 
models such as weakly self-avoiding walks deal with short-range 
potentials.  As far as we know, mathematicians have rarely or never 
studied polyelectrolytes, or even models with long-range potentials.  The 
goal of this paper is to give a mathematically rigorous result related to 
one of the assertions of de Gennes, concerning the spread or radius of 
gyration of the polyelectrolyte.  He uses physical arguments which are not 
mathematically rigorous. We restrict ourselves to the physical case of 
three dimensions, and work in continuous time.  

Now we give the details of our model.  
For $T>0$, let $(B_t)_{0\le t\le T}$ be a standard Brownian motion taking 
values in $\mathbb{R}^3$ with a corresponding probability space 
$(\Omega_T,\mathcal{F}_T,P_T)$.  We denote the usual filtration as 
$(\mathcal{F}_T^{(t)})_{t\in[0,T]}$.  In this paper, we only deal with 
a three-dimensional Brownian motion.  Furthermore, we assume that 
$\Omega_T$ is the canonical probability space, namely the set of 
continuous functions $\omega:[0,T]\to\mathbb{R}^3$ with $\omega_0=0$.  
In that case, $B_t(\omega)=\omega_t$.  

For $\beta>0$, we define a penalization term $\mathcal{E}_T$ and a 
partition function $Z_T$ as follows:
\begin{equation*}
\begin{split}
\mathcal{E}_T = \mathcal{E}_{\beta,T} 
  &:= \exp\left(-\beta\int_{[0,T]^2}\frac{1}{|B_t-B_s|}\dif s\dif t\right)  \\
Z_T = Z_{\beta,T} &:= E^{P_T}\left[\mathcal{E}_{\beta,T}\right].  
\end{split}
\end{equation*}

We will use the convention that $E^{P_T}$ denotes the expectation with 
respect to the probability $P_T$, and likewise for other probabilities.  
Then we define the penalized probability $Q_T$ for events $A\in\mathcal{F}_T$ 
as 
\begin{equation*}
Q_T(A) = Q_{\beta,T}(A) 
  := \frac{1}{Z_{\beta,T}}E^{P_T}\left[\mathcal{E}_{\beta,T}\mathbf{1}_A\right].
\end{equation*}
Of course, $Q_T$ implicitly depends on $\beta$.

Next, we define the radius $R_T$ as
\begin{equation} \label{eq:radius}
R_T := \left(\frac{1}{T^2}\int_{[0,T]^2}|B_t-B_s|^2\dif s\dif t\right)^{1/2}.
\end{equation}

For discrete polymers in three dimensions, de Gennes states that the 
polymer is ``fully extended'', see \cite{deGennes79}, Section XI.2.1, page 299.  
Motivated by this assertion, we prove the following theorem.  

\begin{theorem} \label{th:main}
For positive constants $\beta$ and $c<(2/3)\cdot7^{-2/3}$, we have
\begin{equation}\label{eq:main thm}
\lim_{T\to\infty}Q_T\bigg(c\beta^{1/3}T(\log T)^{-2/3}
   \le R_T\le 6\cdot 7^{1/3}\beta^{1/3}T(\log T)^{1/3}\bigg) = 1.
\end{equation}
\end{theorem}
Inspired by Bolthausen \cite{Bolthausen90} and some previous papers \cite{MN22,MN25}, 
we use the following ideas to prove Theorem \ref{th:main}. Let $\lambda>0$ be a constant,
and define
\begin{equation*}
A_T^{(<)}(\lambda) := \{R_T < \lambda\}
\end{equation*}
and
\begin{equation*}
A_T^{(>)}(\lambda) := \{R_T > \lambda\}.
\end{equation*}
For convenience, we write
\begin{equation*}
\begin{split}
p^{(<)}_T(\lambda) &:= E^{P_T}\Big[\mathcal{E}_T\mathbf{1}_{A_T^{(<)}(\lambda)}\Big],  \\
p^{(>)}_T(\lambda) &:= E^{P_T}\Big[\mathcal{E}_T\mathbf{1}_{A_T^{(>)}(\lambda)}\Big].  
\end{split}
\end{equation*}
Considering the definition of $Q_T$, our goal is to
\begin{enumerate}
\item Bound $p^{(<)}_T(\lambda)$ from above.  
\item Bound $p^{(>)}_T(\lambda)$ from above.  
\item Bound $Z_T$ from below.  
\end{enumerate}

Our proof is short and relies on simple ideas.  In addition, techniques 
like this typically work only in one dimension, although the lace 
expansion (see \cite{MS13}, Chapter 5) works in high dimensions.  
\section{Proof of Theorem \ref{th:main}}

\subsection{The upper bound on $p^{(<)}_T(\lambda)$}

We use H\"older's inequality as follows.  Suppose $a>0$ and $p,q>1$ are 
conjugate exponents, so that $\frac{1}{p}+\frac{1}{q}=1$.  Also, suppose 
that the event $A^{(<)}_T(\lambda)$ occurs.  Then we have
\begin{equation} \label{eq:Holder}
\begin{split}
1 &= \int_{[0,T]^2}|B_t-B_s|^a\cdot|B_t-B_s|^{-a}\frac{\dif s\dif t}{T^2} \\
&\le \left(\int_{[0,T]^2}|B_t-B_s|^{ap}\frac{\dif s\dif t}{T^2}\right)^{1/p}
 \cdot\left(\int_{[0,T]^2}|B_t-B_s|^{-aq}\frac{\dif s\dif t}{T^2}\right)^{1/q}.  
\end{split}
\end{equation}
Now let
\begin{align*}
a=\frac{2}{3}, && p=3, && q=\frac{3}{2},
\end{align*} 
and raise the terms in \eqref{eq:Holder} to the $q$ power.  We find that
\begin{multline*}
1 \le \left(\int_{[0,T]^2}|B_t-B_s|^2\frac{\dif s\dif t}{T^2}\right)^{1/2}
 \left(\int_{[0,T]^2}\frac{1}{|B_t-B_s|}\frac{\dif s\dif t}{T^2}\right) \\
= R_T \int_{[0,T]^2}\frac{1}{|B_t-B_s|}\frac{\dif s\dif t}{T^2}  
\le \lambda\int_{[0,T]^2}\frac{1}{|B_t-B_s|}\frac{\dif s\dif t}{T^2},  
\end{multline*}
where the final inequality above follows from the assumption that 
$A^{(<)}_T(\lambda)$ occurs, and so $R_T<\lambda$. Thus, assuming that 
$A^{(<)}_T(\lambda)$ occurs, we have
\begin{equation*}
\int_{[0,T]^2}\frac{1}{|B_t-B_s|}\dif s\dif t \ge \lambda^{-1}T^2
\end{equation*}
and therefore 
\begin{equation*}
\mathcal{E}_T \le  \exp\left(-\beta\lambda^{-1}T^2\right).
\end{equation*}
Finally, we get
\begin{equation}\label{eq:est of p<}
p^{(<)}_T(\lambda) = E^{P_T}\Big[\mathcal{E}_T\mathbf{1}_{A^{(<)}_T(\lambda)}\Big]
  \le \exp\left(-\beta\lambda^{-1} T^2\right).  
\end{equation}

\subsection{The upper bound on $p^{(>)}_T(\lambda)$}
Since $\mathcal{E}_T\le1$, we have
\begin{equation*}
p^{(>)}_T(\lambda) \le E^{P_T}\Big[\mathbf{1}_{A_T^{(>)}(\lambda)}\Big] 
   = P_T\Big(R_T>\lambda\Big).  
\end{equation*}

Suppose that $R_T>\lambda$.  Since $R_T$ is the 
square root of the mean square distance between pairs 
$(B_s,B_t)_{s,t\in[0,T]}$, we know 
$|B_t-B_s|>\lambda$ for some $s,t\in[0,T]$. From the 
triangle inequality, we then conclude 
$|B_t|>\lambda/2$ for some $t\in[0,T]$. We can write the Brownian motion as
$B_t=(B^{(1)}_t,B^{(2)}_t,B^{(3)}_t)$. Then we have $|B^{(i)}_t|>\lambda/(2\sqrt{3})$ for some 
$i\in\{1,2,3\}$. The reflection principle for the one-dimensional Brownian 
motion and an elementary estimate of normal probabilities imply that
\begin{equation}\label{eq:est of p>}
\begin{split}
p^{(>)}_T(\lambda) 
&\le 3P_T\Big(\sup_{t\in[0,T]}|B^{(1)}_t|>\lambda/(2\sqrt{3})\Big)\\
&\le 12P_T\big(B^{(1)}_T>\lambda/(2\sqrt{3})\big) \\
&  \le \frac{12}{\sqrt{2\pi}}\cdot\frac{\sqrt{T}}{\lambda/(2\sqrt{3})}\exp\left(-\frac{\big(\lambda/(2\sqrt{3})\big)^2}{2T}\right) \\
&= \frac{24\sqrt{3T}}{\sqrt{2\pi }}\cdot\frac{1}{\lambda}\exp\left(-\frac{\lambda^2}{24T}\right).
\end{split}
\end{equation}

\subsection{The lower bound on $Z_T$}
\label{sec:lower bound}

It is helpful to add a constant drift in the first 
coordinate direction to $B$.  Indeed, under $Q_T$ we expect $|B_t|$ to grow 
linearly with $t$, roughly speaking, consistent with de Gennes' assertion that the polymer should be fully extended.  The choice of the first coordinate direction is arbitrary; any other direction would work equally 
well.  To be precise, for $t\in[0,T]$ let
\begin{equation*}
X_t(\omega) = B_t(\omega) + \kappa t\mathbf{e}_1,
\end{equation*}
where $\mathbf{e}_1$ is the unit vector in the first coordinate direction, 
and $\kappa>0$ is a constant to be specified later. According to 
Girsanov's theorem, $X$ induces a probability $\tilde{P}_T$ on 
$(\Omega_T,\mathcal{F}_T)$ with Radon--Nikodym derivative
\begin{equation*}
\frac{\dif\tilde{P}_T}{\dif P_T}(\omega)
  = \exp\left(\int_{0}^{T}\kappa\dif \omega^{(1)}_t-\frac{1}{2}\int_{0}^{T}\kappa^2\dif t\right)
  = \exp\left(\kappa\omega_T^{(1)}-\frac{\kappa^2T}{2}\right),  
\end{equation*}
where $\omega_t^{(1)}=\omega_t\cdot\mathbf{e}_1$ is the first coordinate of 
$\omega_t$.  (Recall $B_t(\omega)=\omega_t$).  $\dif \tilde{P}_T/\dif P_T$ only 
depends on $\omega^{(1)}$ because there is no drift in any other 
coordinate direction.  

We can express $Z_T$ in terms of $\mathcal{E}_T$ and the above 
Radon--Nikodym derivative as  
\begin{equation*}
Z_T = E^{\tilde{P}_T}\left[\mathcal{E}_T
          \left(\frac{\dif\tilde{P}_T}{\dif P_T}\right)^{-1}\right].
\end{equation*}
Since the natural logarithm is a concave function, Jensen's inequality implies
\begin{equation} \label{eq:logZ}
\begin{split}
\log Z_T &\ge E^{\tilde{P}_T}\left[\log\mathcal{E}_T
      - \log\left(\frac{\dif\tilde{P}_T}{\dif P_T}\right)\right]  \\
&= -\beta E^{\tilde{P}_T}
    \left[\int_{[0,T]^2}\frac{1}{|\omega_t-\omega_s|}\dif s\dif t\right]
  - E^{\tilde{P}_T}\left[\kappa\omega_T^{(1)}-\frac{\kappa^2T}{2}\right] \\
&= -\beta E^{\tilde{P}_T}
       \left[\int_{[0,T]^2}\frac{1}{|\omega_t-\omega_s|}\dif s\dif t\right]
  - \frac{\kappa^2T}{2}.  
\end{split}
\end{equation}
We define and estimate
\begin{equation} \label{eq:def_I_1}
\begin{split}
I(T) &:= E^{\tilde{P}_T}
      \left[\int_{[0,T]^2}\frac{1}{|\omega_t-\omega_s|}\dif s\dif t\right] \\
&= \int_{[0,T]^2}E^{P_T}
  \left[\frac{1}{|B_t-B_s+\kappa(t-s)\mathbf{e}_1|}\right]\dif s\dif t   \\
&\le 2T\int_{[0,T]}E^{P_T}\left[\frac{1}{|B_u+\kappa u\mathbf{e}_1|}\right]\dif u.
\end{split}
\end{equation}
We have used the Markov property of Brownian motion and the change of 
variable $u=t-s$ in the last line above.  Since $B_u+\kappa u\mathbf{e}_1$ 
is a normal random variable with mean $\kappa u\mathbf{e}_1$ and 
covariance matrix $uI_3$, we compute
\begin{equation} \label{eq:B+ue}
E^{P_T}\left[\frac{1}{|B_u+\kappa u\mathbf{e}_1|}\right]
= \int_{\mathbb{R}^3}\frac{1}{|x|}
       \cdot\frac{1}{(2\pi u)^{3/2}}
     \exp\bigg(-\frac{|x-\kappa u\mathbf{e}_1|^2}{2u}\bigg)\dif x.
\end{equation}
Note that
\begin{equation} \label{eq:1_over_x}
    \frac{1}{|x|} = \frac{1}{\sqrt{\pi}}\int_0^\infty s^{-1/2}e^{-s|x|^2}\dif s.
\end{equation}
Combining \eqref{eq:B+ue} and \eqref{eq:1_over_x}, we find
\begin{equation} \label{eq:int_1}
\begin{split}
    &E^{P_T}\left[\frac{1}{|B_u+\kappa u\mathbf{e}_1|}\right]
    = \int_{\mathbb{R}^3}\frac{1}{|x|}
       \cdot\frac{1}{(2\pi u)^{3/2}}
       \exp\bigg({-\frac{|x-\kappa u\mathbf{e}_1|^2}{2u}}\bigg)\dif x\\
    &\phantom{E}= \int_{\mathbb{R}^3}\bigg(\frac{1}{\sqrt{\pi}}\int_0^\infty s^{-1/2}e^{-s|x|^2}\dif s\bigg)
       \frac{1}{(2\pi u)^{3/2}}\exp\bigg({-\frac{|x-\kappa u\mathbf{e}_1|^2}{2u}}\bigg)\dif x\\
    &\phantom{E}=\frac{1}{\sqrt{\pi}}\int_0^\infty s^{-1/2}
      \bigg[\frac{1}{(2\pi u)^{3/2}}\int_{\mathbb{R}^3}
      \exp\bigg({-s|x|^2-\frac{|x-\kappa u\mathbf{e}_1|^2}{2u}}\bigg)
      \dif x\bigg]\dif s.
\end{split}
\end{equation}
To continue, we simplify the exponent in the last line of \eqref{eq:int_1} 
by completing the square.  
\begin{equation*}
\begin{split}
-s|x|^2-\frac{|x-\kappa u\mathbf{e}_1|^2}{2u}
&= -\frac{1}{2u}\Big[2us|x|^2+|x|^2-2\kappa u(x\cdot\mathbf{e}_1)+\kappa^2u^2\Big]  \\
&= -\frac{1+2us}{2u}\,\Big|x-\frac{\kappa u}{1+2us}\mathbf{e}_1\Big|^2 \\
&\qquad  + \frac{1+2us}{2u}\left(\frac{\kappa^2 u^2}{(1+2us)^2}-\frac{\kappa^2u^2}{1+2us}\right)  \\
&=:  A_1(x,u,s) + A_2(u,s),
\end{split}
\end{equation*}
where we can further simplify
\begin{equation*}
A_2(u,s) = -\,\frac{\kappa^2su^2}{1+2us}.  
\end{equation*}
Using the fact that
\begin{equation*}
\int_{\mathbb{R}^3}\exp\left(A_1(x,u,s)\right)\dif x 
= \left(\frac{2\pi u}{1+2us}\right)^{3/2},
\end{equation*}
we see that
\begin{equation*}
\begin{split}
E^{P_T}&\left[\frac{1}{|B_u+\kappa u\mathbf{e}_1|}\right]
     = \frac{1}{\sqrt{\pi}}\int_0^\infty s^{-1/2}  \\
&\times \bigg[\frac{1}{(2\pi u)^{3/2}}
    \int_{\mathbb{R}^3}\exp\bigg(A_1(x,u,s)\bigg)\dif x\bigg]
       \exp\big(A_2(u,s)\big)\dif s \\
&\hspace{1cm}=\frac{1}{\sqrt{\pi}}\int_0^\infty s^{-1/2}(1+2us)^{-3/2}
     \exp\bigg(-\frac{\kappa^2su^2}{1+2us}\bigg)\dif s.
\end{split}
\end{equation*}
Making the change of variable $t^2 = \kappa^2su^2/(1+2us)$ and noting that
\begin{equation*}
\dif t = \frac{\kappa u}{2}s^{-1/2}(1+2us)^{-3/2}\dif s,
\end{equation*}
we get
\begin{equation*}
    E^{P_T}\left[\frac{1}{|B_u+\kappa u\mathbf{e}_1|}\right] 
    = \frac{2}{\kappa u\sqrt{\pi}}\int_0^{\kappa\sqrt{u/2}}e^{-t^2}\dif t.
\end{equation*}
Using $\int_{0}^{\kappa\sqrt{u/2}}e^{-t^2}\dif t\le\int_{0}^{\kappa\sqrt{u/2}}\dif t=\kappa\sqrt{u/2}$
and $\int_0^\infty e^{-t^2}\dif t=\sqrt{\pi}/2$, we conclude that
\begin{equation} \label{eq:min}
    E^{P_T}\left[\frac{1}{|B_u+\kappa u\mathbf{e}_1|}\right] 
\le \sqrt{\frac{2}{\pi u}}\mathbf{1}_{[0,\pi/2\kappa^2]}(u)
  + \frac{1}{\kappa u}\mathbf{1}_{[\pi/2\kappa^2,\infty]}(u).  
\end{equation}
Now \eqref{eq:def_I_1} and \eqref{eq:min} imply
\begin{equation} \label{IT-est}
  \begin{split}
  I(T) &\le 2T\int_{[0,T]}E^{P_T}\left[
    \frac{1}{|B_u +\kappa u\mathbf{e}_1|}\right]\dif u  \\
  &\le 2T\int_{0}^{(\pi/2\kappa^2)\wedge T}\sqrt{\frac{2}{\pi u}}\dif u
    +2T\int_{(\pi/2\kappa^2)\wedge T}^{T}\frac{1}{\kappa u}\dif u.  
  \end{split}
\end{equation}
We continue \eqref{IT-est} below.  

In the following, assume that 
\begin{equation} \label{eq:kappa-def}
\kappa=(7\beta\log T)^{1/3}.  
\end{equation}
Note that there exists $T_0>0$ such that $T>T_0$ implies that 
\begin{enumerate}
\item $\kappa>\sqrt{\pi/2T}$, and
\item $4\log \kappa<\log T$.  
\end{enumerate}
Thus, for $T>T_0$, we have $\kappa>\sqrt{\pi/2T}$  
and therefore, $(\pi/2\kappa^2)\wedge T=\pi/2\kappa^2$.  Then,
using \eqref{IT-est}, we have
\begin{equation} \label{eq:IT-est-more}
\begin{split}
I(T) &= 4\kappa^{-1}T+2\kappa^{-1}T\Big[\log T 
        - \log([\pi/2\kappa^2]\wedge T)\Big]  \\
  &\le 6\kappa^{-1}T\log T + 4\kappa^{-1}T\log\kappa \\  
  &\le 7\kappa^{-1}T\log T.  
\end{split}
\end{equation}
Here we have dropped the negative term $-4\log(\pi/2)$ to get an upper 
bound.  

Using \eqref{eq:logZ}, \eqref{eq:def_I_1}, and \eqref{eq:IT-est-more}, we 
find
\begin{equation} \label{eq:lower-bound-log-Z}
\log Z_T \ge 
 -7\beta\kappa^{-1}T\log T -\frac{\kappa^2T}{2}.  
\end{equation}
In fact, we chose the value of $\kappa$ in \eqref{eq:kappa-def} because it 
maximizes the right side of \eqref{eq:lower-bound-log-Z}.  
Substituting $\kappa=(7\beta\log T)^{1/3}$ into 
\eqref{eq:lower-bound-log-Z}, we find that for $T>T_0$,
\begin{equation}\label{eq:est of Z_T}
    Z_T\ge \exp\left(-(3/2)(7\beta)^{2/3}T(\log T)^{2/3}\right).
\end{equation}

\subsection{Completion of the proof of Theorem \ref{th:main}}

Let 
\begin{equation*}
\lambda_1 = c\beta^{1/3}T(\log T)^{-2/3},
\end{equation*}
where $0<c<(2/3)\cdot7^{-2/3}$, and
\begin{equation*}
\lambda_2 = 6\cdot 7^{1/3}\beta^{1/3}T(\log T)^{1/3}.
\end{equation*}
To prove \eqref{eq:main thm}, it suffices to show that
\begin{equation}\label{eq:est of complement}
    \lim_{T\to\infty}Q_T\Big(A^{(<)}_T(\lambda_1)\Big) = \lim_{T\to\infty}Q_T\Big(A^{(>)}_T(\lambda_2)\Big) = 0.
\end{equation}
By \eqref{eq:est of p<} and \eqref{eq:est of Z_T}, for $T>T_0$,
\begin{equation*}
\begin{split}
Q_T\Big(A^{(<)}_T(\lambda_1)\Big) 
&= \frac{E^{P_T}\Big[\mathcal{E}_T\mathbf{1}_{A_T^{(<)}(\lambda_1)}\Big]}{Z_T} = \frac{p_T^{(<)}(\lambda_1)}{Z_T}\\ 
&\le \exp\Big((3/2)7^{2/3}\beta^{2/3}T(\log T)^{2/3}-\beta\lambda_1^{-1}T^2\Big)\\
&= \exp\Big(\big[(3/2)7^{2/3}-c^{-1}\big]\beta^{2/3}T(\log T)^{2/3}\Big)\\
&\to 0 \text{ as $T\to\infty$}.
\end{split}
\end{equation*}
Now by \eqref{eq:est of p>},
\begin{multline*}
Q_T\Big(A^{(>)}_T(\lambda_2)\Big) = \frac{E^{P_T}\Big[\mathcal{E}_T\mathbf{1}_{A_T^{(>)}}(\lambda_2)\Big]}{Z_T}
   = \frac{p_T^{(>)}(\lambda_2)}{Z_T}\\
\le O\left(\lambda_2^{-1}\sqrt{T}\right)
  \exp\bigg((3/2)7^{2/3}\beta^{2/3}T(\log T)^{2/3}-\frac{\lambda_2^{2}}{24T}\bigg)\\
\le O\left(T^{-1/2}(\log T)^{-1/3}\right).
\end{multline*}
Because the exponent in the preceding display is 0.

Therefore, we finally deduce 
\eqref{eq:est of complement} and hence Theorem \ref{th:main}.  

\par\bigskip\noindent
{\bf Acknowledgment.}
We would like to thank Haotian Gu for suggesting an 
improvement to the lower bound in Subsection 
\ref{sec:lower bound}.

\bibliographystyle{amsplain}

\providecommand{\bysame}{\leavevmode\hbox to3em{\hrulefill}\thinspace}
\providecommand{\MR}{\relax\ifhmode\unskip\space\fi MR }
\providecommand{\MRhref}[2]{%
  \href{http://www.ams.org/mathscinet-getitem?mr=#1}{#2}
}
\providecommand{\href}[2]{#2}

\end{document}